\documentclass[twocolumn]{article}
\usepackage[margin=1in]{geometry}
\usepackage{amssymb,graphicx,url,amsmath
}
\usepackage{times}
\usepackage{tikz}
\usepackage{soul}
\usepackage{color}
\usepackage{xcolor}
\usepackage{authblk}

\usepackage{amsmath}
\usepackage{amsfonts}
\usepackage{graphicx}
\usepackage{tikz}
\usepackage{hyperref}
\usepackage{relsize}
\usepackage{breqn}
\usepackage{tcolorbox}
\tcbuselibrary{theorems}
\usetikzlibrary{positioning}
\usetikzlibrary{arrows.meta}
\usetikzlibrary{calc}

\usepackage{verbatim}   
\usepackage{lipsum}

\newif\ifshowtikz
\showtikztrue

\let\oldtikzpicture\tikzpicture
\let\oldendtikzpicture\endtikzpicture

\renewenvironment{tikzpicture}{%
    \ifshowtikz\expandafter\oldtikzpicture%
    \else\comment%
    \fi
}{%
    \ifshowtikz\oldendtikzpicture%
    \else\endcomment%
    \fi
}

\usepackage{xcolor}

\newsavebox\MBox

\newcommand{\bX}{{\bf X}}
\newcommand{\bx}{{\bf x}}

\newcommand{\E}{\mathbb{E}}
\newcommand{\cov}{{\rm Cov}}

\newcommand{\V}{\mathbb{V}}
\newcommand{\Z}{\mathbb{Z}}
\newcommand{\R}{\mathbb{R}}

\newtheorem{Theorem}{Theorem}

\title{Structural causal models for macro-variables in time-series}

\author[1]{Dominik Janzing}
\author[1,2]{Paul Rubenstein} 
\author[1]{Bernhard Sch\"olkopf}
\affil[1]{Max Planck Institute for Intelligent Systemes, T\"ubingen, Germany}
\affil[2]{Machine Learning Group, University of Cambridge, UK}

\date{April 11, 2018}

\begin{document}

\maketitle

\begin{abstract}
We consider a bivariate time series $(X_t,Y_t)$ that is given by a simple linear autoregressive model. Assuming that the equations describing
each variable as a linear combination of past values are considered structural equations, there is a clear meaning of how intervening on
one particular $X_t$ influences 
$Y_{t'}$ at later times $t'>t$. In the present work, we describe conditions under which one can define
a causal model between variables that are coarse-grained in time, thus admitting statements like `setting $X$ to $x$ changes $Y$ in a certain way' without referring to specific time instances. We show that particularly
simple statements follow in the frequency domain, thus providing meaning to interventions on frequencies.
\end{abstract}         
         
\section{Structural equations from dynamical laws}         
         
Structural equations, also called `functional causal models' \cite{Pearl:00}
are a popular and helpful formalization of causal relations. 
For a causal directed acyclic graph (DAG) with $n$ random variables 
$X_1,\dots,X_n$ as nodes they read
\begin{equation}\label{eq:se}
X_j = f_j(PA_j,N_j),
\end{equation}
where $PA_j$ denotes the vector of all parent variables and $N_1,\dots,N_n$ are jointly independent noise variables. Provided the variables $X_j$ refer to measurements that are well-localized in time and correspond to time instances $t_1,\dots,t_n$, one then assumes that $PA_j$ contain only those
variables $X_i$ for which $t_i<t_j$.\footnote{Einstein's special theory of
relativity implies even stronger constraints: If the measurements are also localized in space and correspond to spatial coordinates $z_j$,  
 then $(t_j -t_i)c \geq \|z_j-z_i\|$ where $c$ denotes the speed of light.
 That is, $X_j$ needs to be contained in the forward light cone of all $PA_j$.}           
However, thinking of random variables as measurements that refer to well-defined time instances is too restrictive for many purposes.        
Random variables may, for instance, describe values attained by a quantity
when some system is in its equilibrium state
\cite{Dash05,DGL}. In that case, intervening on one quantity may 
change the stationary joint state, and thus also change the values of other quantities. 

The authors of \cite{DGL} show how the differential equations describing the 
dynamics of the system entail, under fairly restrictive conditions, structural equations relating observables in equilibrium. 
It should be noted, however, that these structural equations may contain causal cycles \cite{Spirtes95,Koster96,PearlDechter96,VoortmanDashDruzdzel10,NIPS_cyclic,Hyttinen12}, i.e., unlike the framework in \cite{Pearl:00} they do not correspond to a DAG. 

The work \cite{Rubenstein16} generalized \cite{DGL}, assaying whether the SCM framework can be extended to model systems that do not converge to an equilibrium (cf.\ also \cite{VoortmanDashDruzdzel10}), and what assumptions need to be made on the ODE and interventions so that this is possible.

Further, also Granger causality \cite{Granger1969} yields coarse-grained statements on causality (subject to appropriate assumptions such as causal sufficiency of the time series) by stating that $X$ causes $Y$ without reference to specific time instances.

The authors of \cite{Chalupka16} propose an approach for the identification of macro-level causes and effects from high-dimensional micro-level measurements
in a scenario that does not refer to time-series. 

In the present work, we will elaborate on the following question: suppose we are given 
 a dynamical system that has a clear causal
 direction when viewed on its elementary time scale. Under which conditions does it
also admit a causal model on `macro-variables' 
that are obtained by coarse-graining variables
referring to particular time instances?
   
\section{Causal models for equilibrium values --- a negative result \label{sec:negative}}
         
\begin{figure} 
\begin{center}
\begin{tikzpicture}[every node/.style = {circle, minimum size=1.3cm}]
{
  \node[draw](xt) {$X_t$};
  \node[draw,  right of=xt, xshift=1.5cm](yt) {$Y_{t}$};
  \node[draw,  below of=xt, yshift=-1.5cm](xt1) {$X_{t+1}$};
  \node[draw,  below of=yt, yshift=-1.5cm](yt1) {$Y_{t+1}$};
  \node[draw,  below of=xt1, yshift=-1.5cm](xt2) {$X_{t+2}$};
  \node[draw,  below of=yt1, yshift=-1.5cm](yt2) {$Y_{t+2}$};

  \node[draw=none,  above of=xt, yshift=1.5cm](xt0) {$\vdots$};
  \node[draw=none,  above of=yt, yshift=1.5cm](yt0) {$\vdots$};  
  \node[draw=none,  below of=xt2, yshift=-1.5cm](xt3) {$\vdots$};
  \node[draw=none,  below of=yt2, yshift=-1.5cm](yt3) {$\vdots$};  

  \node[draw=none,  xshift=-0.4cm] at ($(xt0)!0.5!(xt)$) {$\alpha$};  
  \node[draw=none,  xshift=-0.4cm] at ($(xt)!0.5!(xt1)$) {$\alpha$}; 
  \node[draw=none,  xshift=-0.4cm] at ($(xt1)!0.5!(xt2)$) {$\alpha$}; 
  \node[draw=none,  xshift=-0.4cm] at ($(xt2)!0.5!(xt3)$) {$\alpha$};   
  \node[draw=none, xshift=-0.3cm, yshift=-0.2cm] at ($(xt0)!0.5!(yt)$) {$\beta$};  
  \node[draw=none, xshift=-0.3cm, yshift=-0.2cm] at ($(xt)!0.5!(yt1)$) {$\beta$}; 
  \node[draw=none,xshift=-0.3cm, yshift=-0.2cm] at ($(xt1)!0.5!(yt2)$) {$\beta$}; 
  \node[draw=none,xshift=-0.3cm, yshift=-0.2cm] at ($(xt2)!0.5!(yt3)$) {$\beta$};   
  \node[draw=none,  xshift=-0.4cm] at ($(yt0)!0.5!(yt)$) {$\gamma$};  
  \node[draw=none,  xshift=-0.4cm] at ($(yt)!0.5!(yt1)$) {$\gamma$}; 
  \node[draw=none,  xshift=-0.4cm] at ($(yt1)!0.5!(yt2)$) {$\gamma$}; 
  \node[draw=none,  xshift=-0.4cm] at ($(yt2)!0.5!(yt3)$) {$\gamma$};   
  
  \draw[-{Latex[length=3mm,width=3mm]}](xt)--(xt1);
  \draw[-{Latex[length=3mm,width=3mm]}](yt)--(yt1);
  \draw[-{Latex[length=3mm,width=3mm]}](xt)--(yt1);      
  \draw[-{Latex[length=3mm,width=3mm]}](xt1)--(xt2);
  \draw[-{Latex[length=3mm,width=3mm]}](yt1)--(yt2);
  \draw[-{Latex[length=3mm,width=3mm]}](xt1)--(yt2);      
  
  \draw[-{Latex[length=3mm,width=3mm]}](xt0)--(xt);
  \draw[-{Latex[length=3mm,width=3mm]}](yt0)--(yt);
  \draw[-{Latex[length=3mm,width=3mm]}](xt0)--(yt);     
  \draw[-{Latex[length=3mm,width=3mm]}](xt2)--(xt3);
  \draw[-{Latex[length=3mm,width=3mm]}](yt2)--(yt3);
  \draw[-{Latex[length=3mm,width=3mm]}](xt2)--(yt3);   
  
}
\end{tikzpicture}
\end{center}
\caption{\label{fig:ar1} AR(1) model where $X$ influences $Y$ but not vice versa.}  
\end{figure}
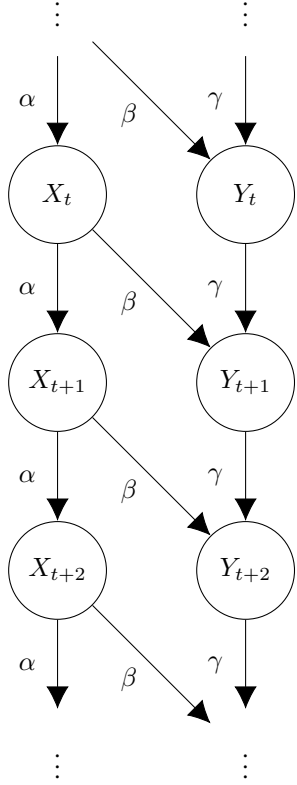

The work \cite{DGL} considered deterministic dynamical systems described by ordinary differential equations and showed that, under particular restrictions, the effect of intervening on
some of the variables changes the equilibrium state of the other ones in a way that can be
expressed by structural equations  among time-less variables,
which are derived from the underlying differential equations. Inspired by these results, we consider non-deterministic discrete
dynamics\footnote{Note that \cite{hansen2014}
considers interventions in stochastic differential equations
and provides conditions under wich they can be seen as limits of interventions
in the `discretized version', such as autoregressive models.}  as given by autoregressive (AR) models,
and ask whether we can define a causal structural equation describing the effect of
an intervention on one variable on another one, which, at the same time, reproduces the
observed stationary joint distribution.
To this end,
we consider the following simple AR model of 
order 1  depicted in Figure~\ref{fig:ar1}.
         
 We assume a Markov chain evolving in time according to the following equations:
\begin{align}
X_{t+1} &= \alpha X_t + E_t^X  \label{eq:sex} \\
Y_{t+1} &= \beta X_t + \gamma Y_t  + E_t^Y 
\label{eq:sey}
\end{align}
Let us assume that $E_t^X,E_t^Y$ are $\mathcal{N}(0,1)$ and \emph{i.i.d.}\ random variables. We assume that the chain goes `back forever' such that $(X_t,Y_t)$ are distributed according to the stationary distribution of the Markov chain, and are jointly normal.\footnote{Note that we could assume initial conditions $X_0$ and $Y_0$, in which case the joint distribution $(X_t,Y_t)$ would not be independent of $t$, but would converge to the stationary distribution.} 

We want to express the stationary distribution and how it changes under (a restricted set of) interventions using a structural causal model. In this example, we consider interventions do($X=x$) and do($Y=y$), by which we refer to the sets of interventions do($X_t=x$) or do($Y_t=y$) for all $t$, respectively.

Let us state our goal more explicitly: we want to derive a structural causal model (SCM) with variables $X$ and $Y$ (and perhaps others) such that the stationary distribution of the Markov chain is the same as the observational distribution on $(X,Y)$ implied by the SCM, and that the stationary distribution of the Markov chain after intervening do($X_t=x$) for all $t$ is the same as the SCM distribution after do($X=x$) (and similar with interventions on $Y$)

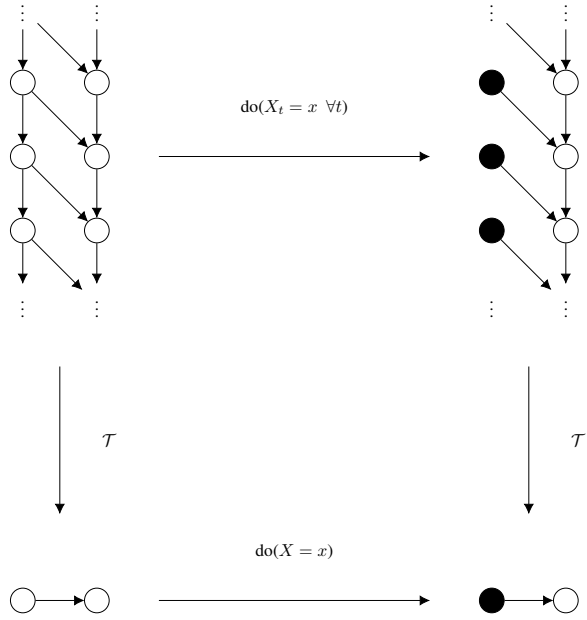
\begin{figure}
\begin{center}
\resizebox{8cm}{!}{
\begin{tikzpicture}[every node/.style = {circle, minimum size=.5cm}]

{
  \node[draw](xt) {};
  \node[draw,  right of=xt, xshift=.5cm](yt) {};
  \node[draw,  below of=xt, yshift=-.5cm](xt1) {};
  \node[draw,  below of=yt, yshift=-.5cm](yt1) {};
  \node[draw,  below of=xt1, yshift=-.5cm](xt2){};
  \node[draw,  below of=yt1, yshift=-.5cm](yt2){};

  \node[draw=none,  above of=xt, yshift=.5cm](xt0) {$\vdots$};
  \node[draw=none,  above of=yt, yshift=.5cm](yt0) {$\vdots$};
  \node[draw=none,  below of=xt2, yshift=-.5cm](xt3) {$\vdots$};
  \node[draw=none,  below of=yt2, yshift=-.5cm](yt3) {$\vdots$};

  \draw[-{Latex[length=2mm,width=2mm]}](xt)--(xt1);
  \draw[-{Latex[length=2mm,width=2mm]}](yt)--(yt1);
  \draw[-{Latex[length=2mm,width=2mm]}](xt)--(yt1);      
  \draw[-{Latex[length=2mm,width=2mm]}](xt1)--(xt2);
  \draw[-{Latex[length=2mm,width=2mm]}](yt1)--(yt2);
  \draw[-{Latex[length=2mm,width=2mm]}](xt1)--(yt2);      
  
  \draw[-{Latex[length=2mm,width=2mm]}](xt0)--(xt);
  \draw[-{Latex[length=2mm,width=2mm]}](yt0)--(yt);
  \draw[-{Latex[length=2mm,width=2mm]}](xt0)--(yt);     
  \draw[-{Latex[length=2mm,width=2mm]}](xt2)--(xt3);
  \draw[-{Latex[length=2mm,width=2mm]}](yt2)--(yt3);
  \draw[-{Latex[length=2mm,width=2mm]}](xt2)--(yt3);   
  
  \node[draw,fill=black,right of=yt, xshift = 7cm](ixt) {};
  \node[draw,  right of=ixt, xshift=.5cm](iyt) {};
  \node[draw,fill=black,  below of=ixt, yshift=-.5cm](ixt1) {};
  \node[draw,  below of=iyt, yshift=-.5cm](iyt1) {};
  \node[draw,fill=black,  below of=ixt1, yshift=-.5cm](ixt2){};
  \node[draw,  below of=iyt1, yshift=-.5cm](iyt2){};

  \node[draw=none,  above of=ixt, yshift=.5cm](ixt0) {$\vdots$};
  \node[draw=none,  above of=iyt, yshift=.5cm](iyt0) {$\vdots$};
  \node[draw=none,  below of=ixt2, yshift=-.5cm](ixt3) {$\vdots$};
  \node[draw=none,  below of=iyt2, yshift=-.5cm](iyt3) {$\vdots$};

  \draw[-{Latex[length=2mm,width=2mm]}](iyt)--(iyt1);
  \draw[-{Latex[length=2mm,width=2mm]}](ixt)--(iyt1);      
  \draw[-{Latex[length=2mm,width=2mm]}](iyt1)--(iyt2);
  \draw[-{Latex[length=2mm,width=2mm]}](ixt1)--(iyt2);

  \draw[-{Latex[length=2mm,width=2mm]}](iyt0)--(iyt);
  \draw[-{Latex[length=2mm,width=2mm]}](ixt0)--(iyt);     
  \draw[-{Latex[length=2mm,width=2mm]}](iyt2)--(iyt3);
  \draw[-{Latex[length=2mm,width=2mm]}](ixt2)--(iyt3);

  \node[draw, below of=xt3,yshift=-5cm](X) {};
  \node[draw,  right of=X, xshift=.5cm](Y) {};

  \draw[-{Latex[length=2mm,width=2mm]}](X)--(Y);

  \node[draw,fill=black, below of=ixt3,yshift=-5cm](iX) {};
  \node[draw,  right of=iX, xshift=.5cm](iY) {};

  \draw[-{Latex[length=2mm,width=2mm]}](iX)--(iY);

  \node[draw=none, right of=yt1](Markov-R){};
  \node[draw=none, left of=ixt1](iMarkov-L){};  
  \draw[-{Latex[length=2mm,width=2mm]}](Markov-R)--(iMarkov-L);

  \node[draw=none,yshift=-1cm](Markov-B) at ($(xt3)!0.5!(yt3)$) {};
  \node[draw=none,yshift=1.5cm](SCM-T) at ($(X)!0.5!(Y)$){};
  \draw[-{Latex[length=2mm,width=2mm]}] (Markov-B) -- (SCM-T);
  
  \node[draw=none, right of=Y](SCM-R){};
  \node[draw=none, left of=iX](iSCM-L){};  
  \draw[-{Latex[length=2mm,width=2mm]}](SCM-R)--(iSCM-L);

  \node[draw=none,yshift=-1cm](iMarkov-B) at ($(ixt3)!0.5!(iyt3)$) {};
  \node[draw=none,yshift=1.5cm](iSCM-T) at ($(iX)!0.5!(iY)$){};
  \draw[-{Latex[length=2mm,width=2mm]}] (iMarkov-B) -- (iSCM-T);
  
  \node[draw=none, yshift=1cm] at ($(Markov-R)!0.5!(iMarkov-L)$) {do($X_t=x \enskip \forall t$)};
  \node[draw=none, yshift=1cm] at ($(SCM-R)!0.5!(iSCM-L)$) {do($X=x$)};
  \node[draw=none, xshift=1cm] at ($(Markov-B) !0.5!(SCM-T)$)  {$\mathcal{T}$};
  \node[draw=none, xshift=1cm] at ($(iMarkov-B)!0.5!(iSCM-T)$) {$\mathcal{T}$};
}
\end{tikzpicture}
}
\end{center}
\caption{\label{fig:commute} Visualization of the projection of the time
resolved causal structure onto a time-less presentation. }  
\end{figure}

This is informally represented by the diagram shown in Figure~\ref{fig:commute}. We seek a `transformation' $\mathcal{T}$ of the original Markov chain (itself an SCM) such that interventions on \emph{all} $X_t$ can be represented as an intervention on a single variable, and such that the SCM summarises the stationary distributions. (Note that in fact as we will see, we cannot express this in general without extra variables as confounders.) The diagram should commute, compare also \cite{Rubensteinetal17}.
 
We first compute the stationary joint distribution  of $(X,Y)$. 
Since there is no influence of $Y$ on $X$, we can first compute the distribution of $X$ regardless of its causal link to $Y$.
Using  
\begin{equation}
 X_{t+1} = E_t^X + \alpha E_{t-1}^X + \alpha^2 E_{t-2}^X + \ldots 
= \sum_{k=0}^\infty \alpha^k E_{t-k}^X, 
\end{equation}  
 and the above conventions
\[
 \E [E^X_t] =0  \quad \hbox{ and } \quad \V [E^X_t] =1,
\] 
we then obtain 
\[
\E [ X_t ] =0    
\]
and
\begin{eqnarray*}
\V[X_t] &=& \mathbb{V}\left[\sum_{k=0}^\infty \alpha^k E_{t-k}^X\right] \\
&=& \sum_{k=0}^\infty \alpha^{2k} \mathbb{V}\left[ E_{t-k}^X\right] \\
&=& \sum_{k=0}^\infty \alpha^{2k}\\
&=& \frac{1}{1-\alpha},
\end{eqnarray*}
where we have used the independence of the noise terms for different $t$.

For the expectation of $Y_t$ we get
\[
\E [Y_t] = \beta \E [X_t ] + \gamma \E [Y_t] + \E [E_t^X] = 0.
\]
To compute the variance of $Y_t$ we need to sum the variance of all independent noise variables. We obtain (the calculation can be found in the appendix):
\[
\V[ Y_t] =\frac{1}{1-\gamma^2} + \frac{\beta^2
(1 + \alpha\gamma)}
{(1-\alpha^2)(1-\alpha\gamma)(1-\gamma^2)} 
\]
For the covariance of $X$ and $Y$ we get (see also appendix):
\[
\cov \left[X_{t},Y_{t}\right] = \frac{\beta\alpha}{(1-\alpha\gamma)(1-\alpha^2)}.
\]
We have thus shown that the stationary joint distribution of $X,Y$ is
\begin{equation}
(X,Y) \sim \mathcal{N}(0,C),
\end{equation}
where the entries of $C$ read
\begin{align*}
C_{XX} &= \frac{1}{1-\alpha^2}\\
C_{XY} &= \frac{\alpha \beta}{(1-\alpha\gamma)(1-\alpha^2)}\\
C_{YY} &= \frac{1}{1-\gamma^2}\\
& + \frac{\beta^2}{(\alpha - \gamma)^2} \left[ \frac{\alpha^2}{1-\alpha^2} - \frac{2\alpha\gamma}{1-\alpha\gamma} + \frac{\gamma^2}{1-\gamma^2}  \right].
\end{align*}
Since the DAG in Figure~\ref{fig:ar1} contains arrows from $X$ to $Y$ and none in the opposite direction,
one would like to explain this bivariate joint distribution by the causal DAG in Figure~\ref{fig:commute} (bottom) where $X$ is causing $Y$. 
This would imply $P(Y|do(X))=P(Y|X)$.
The conditional $P(Y|X)$ is given by a simple regression which yields
\[
Y = a X + E_Y,
\] 
where $E_Y$ is an independent Gaussian noise variable and $a$ is the regression coefficient defined by
\begin{equation}\label{eq:a}
a := C_{XY} C_{XX}^{-1} = \frac{\alpha \beta}{1-\alpha \gamma}. 
\end{equation}
We now show that \eqref{eq:a} is not consistent
with the effect of interventions on $X$ when the latter are defined by setting all $X_t$ to some value $x$.  We refer to this intervention 
as $do(X=x)$. The corresponding interventional distribution of $Y$ reads:
\begin{align*}
Y^{\text{do}(X=x)}_{t+1} &= \beta x + \gamma Y^{\text{do}(X=x)}_{t} + E_t^Y \\
&= \beta x + \beta \gamma x + \beta \gamma^2 x + \ldots +\\
& \quad  E_t^Y + \gamma E_{t-1}^Y + \gamma^2 E_{t-2}^Y + \ldots
\end{align*}
If the distribution is stationary, we have
\[
Y^{\text{do}(X=x)}_{t} = \frac{\beta x}{1-\gamma} + \sum_{k=0}^{\infty} \gamma^k E_{t-k}^Y. 
\]
Hence,
\begin{align*}
Y^{\text{do}(X=x)}_{t} &\sim \mathcal{N}\left(\frac{\beta x}{1-\gamma}, \frac{1}{1-\gamma^2}\right).
\end{align*}
Note that this interventional conditional requires a structural equation whose regression coefficient reads
\begin{equation}\label{eq:a'}
a' := \frac{\beta x}{1-\gamma}, 
\end{equation}
which does not coincide with the coefficient $a$ given by \eqref{eq:a}. 


We now want to provide an intuition about the mismatch between
the regression coefficient $a$  that would be needed
to explain the observed stationary joint distribution and the coefficient $a'$ describing the true effect of interventions.
One the one hand, it should be noted that the conditional of $Y$ given $X$ in the stationary distribution refers to observing only the current
value $X_t$. More precisely, $a$ describes the conditional $P(Y_t|X_t)$, that is, how the distribution of the current value $Y_t$
changes after the current value $X_t$ is observed.
In contrast, the interventions we have considered above are of the type {\em set all variables $X_{t'}$ with $t'\in \mathbb{Z}$ to some fixed value $x$}. 
In other words, the intervention is not localized in time while the observation refers to one specific time instance. 

This motivates already the idea of the following section: in order to explain the observational stationary joint distribution
by an arrow from $X$ to $Y$, we need to define
variables that are de-localized in time because
in that case observations and interventions are de-localized in time.

\section{Non-local variables and frequency analysis}

To understand the reason for the negative result of the preceding section, we 
recall that we compared the interventional conditional variable $Y_t^{do(X=x)}$ (where the intervention $do(X=x)$ referred to all variables $X_t$) to the observational conditional
$Y|X_t=x_t$ (where the observation referred only to the current value $x_t$). 
To overcome this mismatch of completely
non-local interventions on the one hand and
entirely local observations on the other hand, we need to use non-local variables for observations and interventions. This motivates the
following.
For any functions $f,g \in l^1(\Z)$ we define the random variables\footnote{Since $\E[|X_t|]$ and $\E[|Y_t|]$ exist, the series converge in $L^1$ norm of the underlying probability space, hence they converge in probability by Markov's inequality.} 
\[
X_f:= \sum_{t\in \Z} f(t) X_t \quad \hbox{ and } \quad Y_g:=  \sum_{t\in \Z} g(t) Y_t.
\]
One may think of $f,g$ as {\it smoothing kernels} (for instance, discretized Gaussians). Then $X_f,Y_g$ may be the resultand the above conventions
of measurements that perform coarse-graining in time. 
Alternatively, one could also think of $f,g$ as trigonometric functions $\sin,\cos$
restricted to a certain time window. Then
$X_f,Y_g$ are Fourier transforms of the
observations in the respective window. In the spirit of
\cite{Rubensteinetal17}, $X_f$ and $Y_g$ can be thought of as macro-variables derived from the 
micro-variables $X_t,Y_t$ by a `projection'.  
We will now show how a simple causal model
emerges between the macro-variables provided that we consider {\it appropriate pairs}
of macro-variables $X_f,Y_g$.

First, we also define averages over noise variables, which we think of as `macro-noise-variables':
\[
E^X_f:= \sum_{t\in \Z} f(t) E^X_t \quad \hbox{ and } \quad E^Y_g:=  \sum_{t\in \Z} g(t) E^Y_t.
\]
Introducing the shift operator $S$ by $(Sf)(t):=f(t+1)$ we can rewrite 
\eqref{eq:sex} and \eqref{eq:sey}
concisely as
\begin{eqnarray}\label{eq:sexshift}
X_f &=& X_{\alpha S f} + E^X_f \\ \label{eq:seyshift}
Y_g &=& X_{\gamma S g} + Y_{\beta S g} + E_g,
\end{eqnarray}
which can be transformed to
\begin{eqnarray}
X_{(I- \alpha S) f} &=& E^X_f\\
Y_{(I-\beta S)g} &=& X_{\gamma Sg} + E^Y_g,
\end{eqnarray}
and, finally,
\begin{eqnarray}\label{eq:seXf}
X_f &=& E^X_{(I-\alpha S)^{-1} f}\\
\label{eq:seYg}
Y_g &=& X_{\gamma S (I-\beta S)^{-1} g} + E^Y_{(I-\beta S)^{-1} g}. 
\end{eqnarray}
Note that the inverses can be computed from the formal von Neumann series
\[
(I-\alpha S)^{-1} = \sum_{j=0}^\infty (\alpha S)^j,
\]
and $\sum_{j=1}^\infty (\alpha S)^j f$
converges in $l^1(\Z)$-norm
for $\alpha<1$ due to
$\|S^jf\|_1 = \|f\|_1$, and likewise for $\beta<1$. Equation \eqref{eq:seXf} 
describes how the scalar quantity $X_f$
is generated from a single scalar noise term that, in turn, is derived from a weighted average over local noise terms.
Equation \eqref{eq:seYg} describes how the
scalar quantity $Y_g$ is generated from
the scalar $X_{\gamma S (I-\beta S)^{-1} g}$
and a scalar noise term. 

\subsection{Making coarse-graining compatible with the causal model}

The following observation is crucial for the right 
choice of pairs of macro-variables: whenever we choose
\begin{equation}\label{eq:fgrel}
f_g:= \gamma S (I-\beta S)^{-1} g,
\end{equation}
equations~\eqref{eq:seXf} and \eqref{eq:seYg}
turn into the simple form
\begin{eqnarray}\label{eq:seXfsimple}
X_{f_g} &=& E^X_{(I-\alpha S)^{-1} f_g}\\
\label{eq:seYgsimple}
Y_g &=& X_{f_g} + E^Y_{(I-\beta S)^{-1} g}. 
\end{eqnarray}
Equations~\eqref{eq:seXfsimple} and \eqref{eq:seYgsimple} describe how the 
joint distribution of $(X_{f_g},Y_g)$ can be generated: first, generate $X_{f_g}$  from an appropriate average over noise terms. Then, generate $Y_g$ from $X_{f_g}$ plus another averaged noise term. 
For any $x\in \R$, the conditional distribution of $Y_g$, given $X_{f_g}=x$, is therefore identical to the distribution of $x +   E^Y_{(I-\beta S)^{-1} g}$. 

We now argue that \eqref{eq:seXfsimple} and
\eqref{eq:seYgsimple} can even be read as structural equations, that is, they correctly formalize the effect of interventions. 
To this end, we consider the following class
of interventions. For some arbitrary
bounded sequence $\bx= (x_t)_{t\in \Z}$  
we look at the effect of setting ${\bX}$ to ${\bx}$, that is, setting
each $X_t$ to $x_t$. Note that this generalizes
the intervention $do(X=x)$ considered in section~\ref{sec:negative} where each $X_t$ is set to the same value $x\in \R$.  
Using the original structural equation
\eqref{eq:seyshift} yields
\[
Y_g^{do(\bX=\bx)} = \sum_{t\in \Z} x_t f_g(t) +  Y_{\beta S g} +  E_g^Y. 
\]
Applying the same transformations as above yields
\[
Y_g^{do(\bX=\bx)} =  \sum_{t\in \Z} x_t f_g(t)  +  E_{(I-\beta S)^{-1}g}^Y.
\]
Note that the first term on the right hand side
is the value attained by the variable $X_{f_g}$.
Hence, the only information about the entire
intervention that matters for $Y_g$ is the value
of $X_{f_g}$. We can thus talk about `interventions on $X_{f_g}$' without further specifying what the intervention does with each single $X_t$ and write
\[
Y_g^{do(X_{f_g} =x)} = x + E_{(I-\beta S)^{-1}g}^Y.
\]
We have thus shown that \eqref{eq:seXfsimple} and \eqref{eq:seYgsimple}
also reproduce the effect of interventions and
can thus be read as structural equations for the variable pair
$(X_{f_g},Y_g)$.

To be more explicit about the distribution of
the noise terms in \eqref{eq:seXfsimple}
and \eqref{eq:seYgsimple},  
straightforward computation shows the variance to be given by
\begin{align}
\V (E^X_{(I-\alpha S)^{-1}f}) &= \nonumber
\langle (I-\alpha S)^{-1} f,(I-\alpha S)^{-1} f\rangle \\&= \sum_{t\in \Z}    \sum_{k,k'\geq 0} \alpha^{k+k'}  f(t+k') f(t+k). \label{eq:fvar}
\end{align}
Likewise,
\begin{equation}\label{eq:gvar}
\V (E^Y_{(I-\beta S)^{-1}g}) = \sum_{t\in \Z}    \sum_{k,k'\geq 0} \beta^{k+k'}  g(t+k') g(t+k). 
\end{equation}
 We have thus shown the following result.
\begin{Theorem}[valid pairs of macro-variables]
\label{thm:main}
Whenever $f,g\in l^1 (\Z)$ are related by 
\eqref{eq:fgrel}, the AR process in \eqref{eq:sex} and \eqref{eq:sey} entails the scalar structural equations
\begin{eqnarray}
X_f &=& \tilde{E}^X \label{eq:sextheorem}\\
Y_g &=& X_f +  \tilde{E}^Y. \label{eq:seytheorem}
\end{eqnarray}
Here, $\tilde{E}^X$ and $\tilde{E}^Y$ are
zero mean Gaussians whose variances are given by
\eqref{eq:fvar} and \eqref{eq:gvar},
respectively. 

Equation
\eqref{eq:seytheorem} can be read as  a functional `causal model' or `structural equation' in the sense that it describes both the observational conditional of $Y_g$, given $X_g$ and the
interventional conditional of $Y_g$,  $do(X_f=x)$.
\end{Theorem}
In the terminology of \cite{Rubensteinetal17}, the mapping
from the entire bivariate process $(X_t,Y_t)_{y\in \Z}$ to the macro-variable
pair $(X_f,Y_g)$ thus is an {\em exact transformation}
if $f$ and $g$ are related by \eqref{eq:fgrel}.

\subsection{Revisiting the negative result}

Theorem~\ref{thm:main} provides a simple explanation for our negative result from section~\ref{sec:negative}. To see this, we recall that
we considered the distribution of $Y_t$, which
corresponds to the variable $Y_g$ for
$g = (\dots,0,1,0,\dots)$, where the number $1$
occurs at  some arbitrary position $t$. 
To compute the corresponding $f$ according to
\eqref{eq:fgrel} note that
\[
\gamma S (I-\beta S)^{-1} = 
\gamma S \sum_{j=0}^\infty (\beta S)^j =
\frac{\gamma}{\beta} \sum_{j=1}^\infty (\beta S)^j.
\]
We thus obtain the following `smoothing function' $f$ that defines the appropriate coarse graining for $X$ for which
we obtain an {\em exact transformation} of causal models: 
\begin{equation}\label{eq:ffor delta}
f =\gamma (\dots,\beta^2,\beta^1,\beta^0,0,\dots),
\end{equation}
where the first entry from the right is at position $t-1$, in agreement with the intuition that this $X_{t-1}$ is the latest value of $X$ that matters for $Y_t$. 

The intervention $do(X=x)$, where all variables $X_t$ are set to the value $x$, corresponds to
setting $X_f$ to 
\[
x \sum_{t\in \Z} f(t) = x \frac{\gamma}{\beta}.
\]
We thus conclude
\[
Y^{do(X=x)}_t = Y_t^{do(X_f =x \frac{\gamma}{\beta})} = Y_t |_{X_f = x \frac{\gamma}{\beta}}.
\]
In words, to obtain a valid structural equation that formalizes 
both interventional and observational conditional we need to condition
on $X_f$ given by \eqref{eq:ffor delta}.

\subsection{Decoupling of interventions in the frequency domain}

Despite the simple relation between $f$ and $g$
given by \eqref{eq:fgrel},
it is somehow disturbing that {\it different} coarse-grainings are required for $X$ and $Y$.
We will now show that $g$ can be chosen
such that $f$ is {\it almost} the same as $g$
(up to some scalar factor), which
leads us to Fourier analysis of the signals.  

So far, we have thought of $f,g$ as real-valued
functions, but for Fourier analysis it
is instructive to consider complex waves on 
some window $[-T,T]$,
\[
g_{\nu,T} (t):=  \left\{\begin{array}{cc}
 \frac{1}{\sqrt{2T+1}} e^{2\pi i\nu t} &\hbox{ for } t=-T,\dots,T \\
 0  & \hbox{ otherwise.}   
 \end{array} \right. 
\]
For notational convenience, we also introduce the 
corresponding vectors $f$ by
\[
f_{\nu,T} := \gamma S (I-\beta S)^{-1} g_{\nu,T},
\]
which are not as simple as $g_{\nu,T}$.
However, 
for sufficiently large $T$, the functions $g_{\nu,T}$ are 
almost eigenvectors of $S$ with eigenvalue $z_\nu :=e^{2\pi i \nu}$ since we have
\begin{equation}\label{eq:almosteigen}
\|S^j g_{\nu,T} - z_\nu^j g_{\nu,T} \|_1 \leq \frac{2j}{\sqrt{2T+1}},
\end{equation}
because the functions differ only at the positions 
$-T,\dots,-T+j-1$ and $T+1,\dots,T+j$.
We show in the appendix that this implies
\begin{eqnarray}\label{eq:almosteigenderived}
&&\|f_{\nu,T} - \gamma z_\nu (1-z_\nu)^{-1} g_{\nu,T}\|_1\\
&\leq& \nonumber
\frac{2}{\sqrt{2T+1}}  \frac{|\gamma|}{|\beta|}  |(1-\beta)^{-2}|.
\end{eqnarray}
that is, $f_{\nu,T}$ coincides with a complex-valued multiple of $g_{\nu,T}$ up to
an error term that decays with $O(1/\sqrt{T})$. 
Using the abbreviations 
\[
E_{\nu,T}^X := E^X_{(I-\alpha S)^{-1}f_{\nu,T}},
\]
and
\[
E_{\nu,T}^Y := E^Y_{(I-\beta S)^{-1}g_{\nu,T}},
\]
the general structural equations \eqref{eq:seXfsimple} and \eqref{eq:seYgsimple} thus imply 
 the approximate structural 
equations
\begin{align}\label{eq:asx}
X_{g_{\nu,T}} &= E_{\nu,T}^X,\\
Y_{g_{\nu,T}} &\approx  \gamma e^{2\pi i \nu} (1- \beta e^{2\pi i \nu})^{-1} X_{g_{\nu,T}} + E_{\nu,T}^Y,
\label{eq:asy}
\end{align}
where the error of the approximation \eqref{eq:asy} is a random variable
whose $L^1$-norm is bounded by
\[
\frac{2}{\sqrt{2T+1}}  \frac{|\gamma|}{|\beta|}  |(1-\beta)^{-2}| \cdot \E [|X_t|],
\]
due to \eqref{eq:almosteigenderived}. 
We conclude with the interpretation that
the structural equations for different frequencies perfectly {\em decouple}. That is,
intervening on one frequency of $X$ has only effects on the same frequency of $Y$,
as a simple result of linearity and time-invariance of the underlying Markov process. 

To phrase this decoupling over frequencies
in a precise way, show that $ E^X_{\nu,T}$ 
and 
$
  E^Y_{\nu,T} 
$
exist in distribution as complex-valued random variables. It is sufficient to show that
the variances and covariances of real and imaginary parts of $E_{\nu,T}^Y$ converge because both variables are Gaussians with zero mean. We have 
\begin{align}
&\V [{\rm Re} \{E^Y_{\nu,T}\}] = \frac{1}{4}\E [ (E^Y_{\nu,T} +  \bar{E}^Y_{\nu,T})^2] =\\
&
\frac{1}{4}\left( \E [ (E^Y_{\nu,T})^2]  + 
\E [ (\bar{E}^Y_{\nu,T})^2] + 2
\E [ \bar{E}^Y_{\nu,T} E^Y_{\nu,T} ]\right).
\end{align}

We obtain
\begin{align}\label{eq:asymvar}
& \E [ \bar{E}^Y_{\nu,T} E^Y_{\nu,T} ] =   \E[ \overline{E^Y_{(I-\beta S)^{-1} g_{\nu,T}}} E^Y_{(I-\beta S)^{-1}g_{\nu,T}} ]\\
&=
\langle (I+\beta S)^{-1} g_{\nu,T} , 
(I-\beta S)^{-1} g_{\nu,T}\rangle\\
&\to  |(1-\beta z_\nu)^{-1}| ^2. \nonumber  
\end{align}
For the first equality, recall that the complex-valued inner product is anti-linear in its first argument. The limit follows from straightforward computations using an analog of \eqref{eq:almosteigen} for the $L^2$ norm,
\[
\|S^j g_{\nu,T} - z_\nu^j g_{\nu,T} \|_2 \leq \sqrt{\frac{2j}{2T+1}},
\] 
and further algebra akin to the proof of \eqref{eq:almosteigenderived} in the appendix.

Moreover, 
\begin{align}
& \E [ (E^Y_{\nu,T})^2] = 
\E[ E^Y_{(I-\beta S)^{-1} g_{\nu,T}} E^Y_{(I-\beta S)^{-1}g_{\nu,T}} ] =\\
&  \langle \overline{(I+\beta S)^{-1} g_{\nu,T} }, 
(I-\beta S)^{-1} g_{\nu,T}\rangle. 
\end{align} 
Hence, $\E[(E^Y_{\nu,T})^2]$
and its conjugate  $\E[\overline{(E^Y_{\nu,T})^2}]$
converge to zero for all $\nu\neq 0$ because 
\begin{align}
&\lim_{T\to\infty}  \langle \overline{(I+\beta S)^{-1} g_{\nu,T} }, 
(I-\beta S)^{-1} g_{\nu,T}\rangle \label{eq:withS}\\
&(1-\beta z_\nu)^{-2} \lim_{T\to\infty} 
\langle \overline{g_{\nu,T}}, g_{\nu,T} \rangle \label{eq:withz} \\
& =  (1-\beta z_\nu)^{-2} \lim_{T\to\infty} 
\sum_t g^2_{\nu,T} (t) = 0, 
\end{align}
where equality of \eqref{eq:withS} and \eqref{eq:withz} 
follows from \eqref{eq:almosteigen}. 
Hence, only the mixed term containing both $E^Y_{\nu,T}$ and its conjugate survives the limit.
We conclude 
\[
\lim_{T\to \infty} V [{\rm Re} \{E^Y_{\nu,T}\}] 
= \frac{1}{2} |(1-\beta z_\nu)^{-1}|^2. 
\]
Similarly, we can show that  
$V [{\rm Im} \{E^Y_{\nu,T}\}] $ converges to the same value. 
Moreover,
\[
\lim_{T\to\infty} \cov [{\rm Re} \{E^Y_{\nu,T}\} ,{\rm Im} \{E^Y_{\nu,T}\}] =0,
\]
because straightforward computation shows that 
the covariance contains no mixed terms. 
Hence we can define
\[
E_\nu^Y := \lim_{T\to \infty} E^Y_{\nu,T},
\]
with convergence in distribution. Real and imaginary parts are uncorrelated and their variance read: 
\begin{equation}\label{eq:asymvarY}
\V[{\rm Re}\{E^Y_\nu\}] =  \V [{\rm Im} \{E^Y_{\nu}\}]
= \frac{1}{2} |(1-\beta z_\nu)^{-1}|^2.
\end{equation}
We conclude that the distribution of $E^Y_\nu$ is an isotropic Gaussian in the complex plane, whose
components have variance
$\frac{1}{2} |(1-\beta z_\nu)^{-1}|^2$.   

To compute the limit of $E^X_{\nu,T}$ we
proceed similarly and observe
\begin{align}
& \E [ \bar{E}^Y_{\nu,T} E^Y_{\nu,T} ] =   \E[ \overline{E^Y_{\gamma S (I-\beta S)^{-2} g_{\nu,T}}} E^Y_{\gamma S(I-\beta S)^{-2}g_{\nu,T}} ]\nonumber \\
&\to  |\gamma z_\nu (1-\beta z_\nu)^{-2}| ^2.
\end{align} 
We can therefore define the 
random variable $E^X_\nu:= \lim_{T\to \infty}$ (again with convergence in distribution) with 
\begin{equation}\label{eq:asymvarX}
\V[{\rm Re}\{E^X_\nu\}] =  \V [{\rm Im} \{E^X_{\nu}\}]= 
\frac{1}{2} \left|\frac{\gamma z_\nu}{(1-\beta z_\nu)^2}\right|^2.
\end{equation}
We can phrase these findings by asymptotic structural equations
\begin{align*}
X_\nu &=E^X_\nu\\
Y_\nu &= \gamma e^{2\pi i \nu} (1-\beta e^{2\pi i \nu})^{-1} X_\nu + E^Y_\nu,
\end{align*}
where the variances of real and imaginary parts of $E_\nu^X$ and $E_\nu^Y$ are given by \eqref{eq:asymvarX} and
\eqref{eq:asymvarY}, respectively.

\section{Conclusion}
We have studied bivariate time series $(X_t,Y_t)$ given by linear autoregressive models, and described conditions under which one can define
a causal model between variables that are coarse-grained in time, thus admitting statements like `setting $X$ to $x$ changes $Y$ in a certain way' without referring to specific time instances. We show that particularly
elegant statements follow in the frequency domain, thus providing meaning to interventions on frequencies.

\section{Appendix}

\subsection{Covariance of $X$ and $Y$ in the stationary distribution}

\begin{align*}
&\cov \left[X_{t+1},Y_{t+1}\right] \\
&= \mathbb{C}\left[\sum_{k=0}^\infty \alpha^k E_{t-k}^X, \sum_{k=0}^\infty \gamma^k E_{t-k}^Y+ \beta \sum_{k=0}^\infty
 \frac{\alpha^{k+1}-\gamma^{k+1}}{\alpha-\gamma}  E_{t-1-k}^X  \right] \\
 &= \cov \left[\sum_{k=0}^\infty \alpha^{k+1} E_{t-k-1}^X, \beta \sum_{k=0}^\infty
 \frac{\alpha^{k+1}-\gamma^{k+1}}{\alpha-\gamma}  E_{t-1-k}^X  \right] \\
 &= \frac{\beta}{\alpha-\gamma}\sum_{k=0}^\infty \alpha^{2k+2}-\alpha^{k+1}\gamma^{k+1} \\
 &= \frac{\beta}{\alpha-\gamma} \left[ \frac{\alpha^2}{1-\alpha^2} - \frac{\alpha\gamma}{1-\alpha\gamma} \right] \\
 &= \frac{\beta}{\alpha-\gamma} \left[ \frac{\alpha^2(1-\alpha\gamma) - \alpha\gamma(1-\alpha^2)}{(1-\alpha^2)(1-\alpha\gamma)} \right] 
\end{align*}

\subsection{Approximate eigenvalues of functions of $S$}

Using \eqref{eq:almosteigen} we obtain 
\begin{eqnarray*}
&&\left\|\gamma S (I-\beta S)^{-1} g_{\nu,T} - \gamma z_\nu (1-z_\nu)^{-1} g_{\nu,T}\right\|_1\\
 &=&
\left\| \frac{\gamma}{\beta} \sum_{j=1}^\infty (\beta S)^j g_{\nu,T} -
\frac{\gamma}{\beta} \sum_{j=1}^\infty (\beta z_\nu)^j g_{\nu,T}\right\|_1\\   
&\leq & \frac{|\gamma|}{|\beta|}  \left|\sum_{j=1}^\infty \frac{2j \beta^j}{\sqrt{2T+1}}\right|  = \frac{2}{\sqrt{2T+1}}  \frac{|\gamma|}{|\beta|}  \left|\frac{d}{d\beta} \sum_{j=0}^\infty \beta^j\right| \\
&=& \frac{2}{\sqrt{2T+1}}  \frac{|\gamma|}{|\beta|}  \left|\frac{d}{d\beta} (1-\beta)^{-1}\right|  = \frac{2}{\sqrt{2T+1}}  \frac{|\gamma|}{|\beta|}  |(1-\beta)^{-2}|.
\end{eqnarray*}

\subsection{Variance of $Y$ in the stationary distribution}

\begin{align*}
&\mathbb{V}\left[Y_{t+1} \right] \\
&= \mathbb{V}\left[\sum_{k=0}^\infty \gamma^k E_{t-k}^Y \right] + \mathbb{V}\left[ \beta \sum_{k=0}^\infty
 \frac{\alpha^{k+1}-\gamma^{k+1}}{\alpha-\gamma}  E_{t-1-k}^X \right] \\
 &= \sum_{k=0}^\infty \gamma^{2k}+ \frac{\beta^2}{(\alpha-\gamma)^2} \sum_{k=0}^\infty \alpha^{2k+2} -2\alpha^{k+1}\gamma^{k+1}-\gamma^{2k+2}  \\
 &= \frac{1}{1-\gamma^2} + \frac{\beta^2}{(\alpha-\gamma)^2} \left[\frac{\alpha^2}{1-\alpha^2} - \frac{2\alpha\gamma}{1-\alpha\gamma} + \frac{\gamma^2}{1-\gamma^2} \right] \\
& = \frac{1}{1-\gamma^2} + \frac{\beta^2}{(\alpha-\gamma)^2}\times\\
&
\left[ \frac{\alpha^2(1-\alpha\gamma)(1-\gamma^2) -2\alpha\gamma(1-\alpha^2)(1-\gamma^2)}{(1-\alpha^2)(1-\alpha\gamma)(1-\gamma^2)} \right.\\
&
\left.+\frac{\gamma^2(1-\alpha^2)(1-\alpha\gamma)}{(1-\alpha^2)(1-\alpha\gamma)(1-\gamma^2)}\right] \\
& = \frac{1}{1-\gamma^2} + \frac{\beta^2}{(\alpha-\gamma)^2}\times \\
&\left[ \frac{%
\alpha^2 - \textcolor{black}{\alpha^3\gamma} - \textcolor{black}{\alpha^2\gamma^2} + \textcolor{black}{\alpha^3\gamma^3}  - %
2\alpha\gamma + \textcolor{black}{2\alpha^3\gamma}  }{(1-\alpha^2)(1-\alpha\gamma)(1-\gamma^2)}\right.\\ 
&\left.
+\frac{\textcolor{black}{2\alpha\gamma^3} - \textcolor{black}{2\alpha^3\gamma^3} + %
\gamma^2 - \textcolor{black}{\alpha^2\gamma^2} -\textcolor{black}{\alpha\gamma^3} + \textcolor{black}{\alpha^3\gamma^3}}%
{(1-\alpha^2)(1-\alpha\gamma)(1-\gamma^2)}\right] \\
& = \frac{1}{1-\gamma^2} + \frac{\beta^2}{(\alpha-\gamma)^2} \left[ \frac{%
\alpha^2 + \alpha^3\gamma  - %
2\alpha\gamma  + \alpha\gamma^3 + %
\gamma^2 - 2\alpha^2\gamma^2  }%
{(1-\alpha^2)(1-\alpha\gamma)(1-\gamma^2)}\right] \\
& = \frac{1}{1-\gamma^2} + \frac{\beta^2}{(\alpha-\gamma)^2} \left[ \frac{%
\alpha^2 - 2\alpha\gamma + \gamma^2 + \alpha\gamma(\alpha^2 - 2\alpha\gamma    + \gamma^2  }%
{(1-\alpha^2)(1-\alpha\gamma)(1-\gamma^2)}\right] \\
& = \frac{1}{1-\gamma^2} + \frac{\beta^2}{(\alpha-\gamma)^2} \left[ \frac{%
(\alpha-\gamma)^2(1 + \alpha\gamma)}%
{(1-\alpha^2)(1-\alpha\gamma)(1-\gamma^2)}\right] \\
%
%
\end{align*}


\begin{thebibliography}{10}

\bibitem{Pearl:00}
J.~Pearl.
\newblock {\em Causality}.
\newblock Cambridge University Press, 2000.

\bibitem{Dash05}
D.~Dash.
\newblock Restructing dynamic causal systems in equilibrium.
\newblock In {\em Proc. Uncertainty in Artifical Intelligence}, 2005.

\bibitem{DGL}
J.~Mooij, D.~Janzing, and B.~Sch\"olkopf.
\newblock From ordinary differential equations to structural causal models: the
  deterministic case.
\newblock In Nicholson A. and P.~Smyth, editors, {\em Proceedings of the 29th
  Conference on Uncertainty in Artificial Intelligence (UAI)}, pages 440--448,
  Oregon, USA, 2013. AUAI Press Corvallis.

\bibitem{Spirtes95}
P.~Spirtes.
\newblock Directed cyclic graphical representations of feedback models.
\newblock In {\em Proceedings of the Eleventh Conference on Uncertainty in
  Artificial Intelligence}, UAI'95, pages 491--498, San Francisco, CA, USA,
  1995. Morgan Kaufmann Publishers Inc.

\bibitem{Koster96}
J.~T.~A. Koster.
\newblock Markov properties of nonrecursive causal models.
\newblock {\em Annals of Statistics}, 24(5):2148--2177, 1996.

\bibitem{PearlDechter96}
J.~Pearl and R.~Dechter.
\newblock Identifying independence in causal graphs with feedback.
\newblock In {\em Proceedings of the Twelfth Annual Conference on Uncertainty
  in Artificial Intelligence (UAI-96)}, pages 420--426, 1996.

\bibitem{VoortmanDashDruzdzel10}
M.~Voortman, D.~Dash, and M.~Druzdzel.
\newblock Learning why things change: The difference-based causality learner.
\newblock In {\em Proceedings of the Twenty-Sixth Annual Conference on
  Uncertainty in Artificial Intelligence (UAI)}, pages 641--650, Corvallis,
  Oregon, 2010. AUAI Press.

\bibitem{NIPS_cyclic}
J.~Mooij, D.~Janzing, B.~Sch\"olkopf, and T.~Heskes.
\newblock Causal discovery with cyclic additive noise models.
\newblock In {\em Advances in Neural Information Processing Systems 24,
  Twenty-Fifth Annual Conference on Neural Information Processing Systems (NIPS
  2011), Curran}, pages 639--647, NY, USA, 2011. Red Hook.

\bibitem{Hyttinen12}
A.~Hyttinen, F.~Eberhardt, and P.O. Hoyer.
\newblock Learning linear cyclic causal models with latent variables.
\newblock {\em Journal for Machine Learning Research}, 13:3387−3439, November
  2012.

\bibitem{Rubenstein16}
P.~K. Rubenstein, S.~Bongers, J.~M. Mooij, and B.~Sch{\"{o}}lkopf.
\newblock From deterministic {ODEs} to dynamic structural causal models.
\newblock {\em arXiv}, 1608.08028, 2016.

\bibitem{Granger1969}
C.~W.~J. Granger.
\newblock Investigating causal relations by econometric models and
  cross-spectral methods.
\newblock {\em Econometrica}, 37(3):424--38, July 1969.

\bibitem{Chalupka16}
K.~Chalupka, P.~Perona, and F.~Eberhardt.
\newblock Multi-level cause-effect systems.
\newblock In {\em Proceedings of the 19th International Conference on
  Artificial Intelligence and Statistics (AISTATS)}, JMLR: W\&CP volume 41,
  2016.

\bibitem{hansen2014}
N.~Hansen and Al. Sokol.
\newblock Causal interpretation of stochastic differential equations.
\newblock {\em Electron. J. Probab.}, 19:1–24, 2014.

\bibitem{Rubensteinetal17}
P.~K. Rubenstein, S.~Weichwald, S.~Bongers, J.~M. Mooij, D.~Janzing,
  M.~Grosse-Wentrup, and B.~Sch{\"o}lkopf.
\newblock Causal consistency of structural equation models.
\newblock In {\em Proceedings of the Thirty-Third Conference on Uncertainty in
  Artificial Intelligence (UAI)}, 2017.

\end{thebibliography}
\end{document}